\documentclass{article}

\usepackage{latexsym,amsfonts,amsmath,amsthm,graphics}
\usepackage{amssymb,mathtools, url,dsfont}
\usepackage{epsfig}
\usepackage{color}
\usepackage{cellspace}
\usepackage{tikz}

\usepackage{xspace}
\usepackage{enumitem}
\usepackage{pb-diagram}

\usepackage[a4paper,dvips]{geometry}

\usepackage{tocloft}

\usepackage[font=small,labelfont=bf]{caption}

\newtheorem{thm}{Theorem}[section]
\newtheorem{cor}[thm]{Corollary}

\newtheorem{lem}[thm]{Lemma}
\newtheorem{defi}[thm]{Definition}
\newtheorem{prop}[thm]{Proposition}

\newtheorem{rem}[thm]{Remark}
\newtheorem{ex}[thm]{Example}

\newenvironment{Proof}{\begin{trivlist}

\item[\hskip\labelsep{\it Proof.}]}{$\hfill\Box$\end{trivlist}}

\newcommand{\RR}{\mathbb{R}}

\newcommand{\Ss}{\mathbb{S}}

\newcommand{\x}{{\bf x}}
\newcommand{\y}{{\bf y}}
\newcommand{\abf}{{\bf a}}
\newcommand{\bbf}{{\bf b}}

\DeclareMathOperator{\sgn}{sgn}

\newcommand{\R}{\mathbb{R}}

\newcommand{\N}{\mathbb{N}}

\newcommand{\eps}{\varepsilon}

\newcommand{\Y}{\mathcal{Y}}

\newcommand{\Vs}[1] {V_{\mathcal{S},#1 } } 		
\newcommand{\Va}[1] {V_{#1}}						

\newcommand{\Vi}[1] {\mathcal{V}_{\infty}(#1)}	
\newcommand{\Vf}[1] {\mathcal{V}(#1)}			

\begin{document}

\title{On functions of bounded variation}


\author{Christoph Aistleitner\footnote{Johannes Kepler University Linz, Department of Financial Mathematics and Applied Number Theory, Linz, Austria. E-Mail: \emph{aistleitner@math.tugraz.at}}, 
Florian Pausinger\footnote{ IST Austria, Am Campus 1, A-3400 Klosterneuburg, Austria. E-Mail: \emph{florian.pausinger@gmx.at} },
Anne Marie Svane\footnote{Aarhus University, Aarhus, Denmark. E-Mail: amsvane@math.au.dk} and
Robert F. Tichy\footnote{TU Graz, Department for Analysis and Computational Number Theory (Math A), Graz, Austria. E-Mail:\emph{tichy@tugraz.at}}
 }

\maketitle

\begin{abstract}
The recently introduced concept of $\mathcal{D}$-variation unifies previous concepts of variation of multivariate functions. In this paper, we give an affirmative answer to the open question from \cite{PauSva14} whether \emph{every} function of bounded Hardy--Krause variation is Borel measurable and has bounded $\mathcal{D}$-variation. Moreover, we show that the space of functions of bounded $\mathcal{D}$-variation can be turned into a commutative Banach algebra.
\\[3pt]
\\[3pt]
\end{abstract}

\tableofcontents



\section{Introduction}
\label{sec1}

It is a classical problem to generalise the notion of total variation of a one-dimensional function to multivariate functions and  study conditions under which a function has bounded variation.
The algebraic properties of the corresponding spaces of functions of bounded variation are of particular interest in numerical integration.
Let $f$ be a real-valued measureable function over a compact Hausdorff space $X$, equipped with a sigma-field $\mathcal{F}$ and a normalized measure $\mu$. Furthermore, let $x_1, x_2 \ldots, x_N \in X$.
The famous Koksma--Hlawka inequality is a general principle to bound the approximation error 
\begin{equation}
\left| \frac{1}{N} \sum_{j=1}^N f(x_j) - \int_X f(x) ~d\mu(x) \right|
\end{equation}  
by the product of two independent factors. One of these factor depends only on the function $f$ (more precisely, on the \emph{variation} of $f$) and the other factor depends only on the discrete point set (the \emph{discrepancy} of $x_1, x_2, \ldots, x_N$). Informally speaking, the discrepancy measures the deviation between the empirical distribution of the points $x_1, \dots, x_N$ and the measure $\mu$. The classical setting is that of $X$ being the $d$-dimensional unit cube $[0,1]^d$ and $\mu$ being the $d$-dimensional Lebesgue measure; in this setting, the Koksma--Hlawka inequality reads as 
\begin{equation} \label{kh_classic}
\left| \frac{1}{N} \sum_{j=1}^N f(\mathbf{x}_j) - \int_{[0,1]^d} f(\mathbf{x}) ~d\mathbf{x} \right| \leq D_N^*(\mathbf{x}_1, \dots, \mathbf{x}_N) \cdot \textup{Var}_{\textup{HK}} ~f,
\end{equation}
where $D_N^*$ is the so-called star discrepancy and $\textup{Var}_{\textup{HK}}$ is the variation in the sense of Hardy and Krause; for details on this basic result of discrepancy theory, see for example \cite{dts,knu}.\\

It is well known that the space of all real-valued functions of bounded total variation on the compact interval $[a,b]$ is a commutative Banach algebra with respect to pointwise multiplication.
However, it is not obvious how to generalize this notion of bounded variation to the case of multivariate functions. Hardy \cite{hardy} and Krause \cite{krause} introduced a concept of bounded variation for multivariate functions, which was used by Hlawka \cite{hlawka} to generalize the one-dimensional Koksma inequality \cite{koksma} and to obtain the classical version of the Koksma--Hlawka inequality as stated in \eqref{kh_classic}. G\"otz \cite{gotz} proved a version of the Koksma--Hlawka inequality for general measures (rather than only Lebesgue measure), and recently, Brandolini, Colzani, Gigante and Travaglini \cite{brandolini, brandolini2} replaced the integration domain $[0,1]^s$ by an arbitrary bounded Borel subset of $\RR^d$ and proved the inequality for piecewise smooth integrands.\\

The notion of Hardy--Krause variation was generalised in a natural way by Bl\"{u}mlinger and Tichy \cite{BluTic89}, who proved that the corresponding space of functions of bounded variation is a commutative 
Banach algebra.
However, especially in the context of numerical integration, these different notions of Hardy--Krause variation come with the severe drawback that many functions of practical interest have unbounded variation (e.g. the indicator function of a ball or a tilted box). Recently Harman \cite{Har10} introduced a new notion of variation, which remains finite for certain discontinuous functions with unbounded variation in the sense of Hardy and Krause, and proved a Koksma--Hlawka inequality in this settig. Unfortunately, the space of functions of bounded Harman variation lacks many of the nice algebraic properties of the Hardy--Krause variation.\\

This was the motivation that led to the introduction of a general framework of variations in \cite{PauSva14}. The concept of $\mathcal{D}$-variation unifies the different notions of variation and is not restricted to integrals over $[0,1]^d$, but works for integrals over arbitrary compact Hausdorff spaces.
It was shown to coincide with Hardy--Krause variation in special cases.
The first aim of our paper is to show that \emph{every} function of bounded Hardy--Krause variation also has bounded variation in the new sense, thus answering a question which was left open in \cite{PauSva14}; see Section \ref{sec4}.
In particular this also means that every function of bounded Hardy--Krause variation is Borel measurable. This is a fundamental result which we did not find anywhere in the literature, for which reason we also provide a self-contained proof in Section \ref{sec3}.\\

Functions of bounded Hardy--Krause variation received a lot of attention in the literature; see \cite{adams,AisDic15,BluTic89,Blu89,clarkson,leonov}. Given our results, it is natural to ask whether results about the structure of the space of functions of bounded Hardy--Krause variation, e.g., that it is a Banach algebra (see \cite{BluTic89,Blu89}), also extend to our more general notion. We discuss this question in Section \ref{sec5}, where we show that the space of functions of bounded $\mathcal{D}$-variation is indeed a commutative Banach algebra.


\section{Different notions of variation}
\label{sec2}

In the following we introduce the two definitions of variation of a multivariate function that we consider: the classical Hardy--Krause variation and the recently introduced $\mathcal{D}$-variation.

\subsection{Hardy--Krause variation}

{\bf Definition. }
In the following, we use the notation of Owen \cite{Owe05}.
Let $f(\x)$ be a function on $[0,1]^d$. If $\abf=(a_1, \ldots, a_d)$ and $\bbf=(b_1, \ldots, b_d)$ are elements of $[0,1]^d$ such that $a_i \leq b_i \ (a_i < b_i)$ for all $1\leq i \leq d$, then we write $\abf \leq \bbf \ (\abf < \bbf)$. 
For $u\subseteq \{1,\dots,d\}$, we denote by $\abf^u:\bbf^{-u}$ the point with $i$-th coordinate equal to $a_i$ if $i\in u$ and equal to $b_i$ otherwise. The set $-u$ is the set complement of $u$ in $\{1,\dots,d\}$.
Using this notation, we introduce the $d$-dimensional difference operator 
$$\Delta^{(d)} (f;R)=\Delta (f;R) = \sum_{u \subseteq \{1,\ldots, d\}} (-1)^{|u|} f(\abf^u : \bbf^{-u}),$$ 
which assigns to the axis-parallel rectangle $R=[\abf,\bbf]$ a $d$-dimensional quasi-volume.\\

In dimension $d=1$, a ladder $\Y$ on the interval $[0,1]$ is a partition of $[0,1]$, i.e. a sequence $0=y_1 < \dotsm < y_k <1$.
A ladder in $[0,1]^d$ is a set of the form $\mathcal{Y} = \prod_{j=1}^d \mathcal{Y}^j\subseteq [0,1]^d$, where each $\mathcal{Y}^j$ is a one-dimensional ladder. Let $\mathbb{Y}$ be the set of all ladders on $[0,1]^d$.
Suppose $\Y^j=\{y_1^j<\dotsm < y_{k_j}^j\}$. Define the successor $(y_i^j)_+$ of $y_i^j$ to be $y_{i+1}^j$ if $i< k_j$ and $(y_{k_j}^j)_+=1$. If $\y=(y_{i_1}^1,\dots,y_{i_d}^d) \in \mathcal{Y}$, then we define its  successor to be $\y_+= ((y_{i_1}^1)_+,\dots,(y_{i_d}^d)_+)$.
For a ladder $\mathcal{Y}$ in $[0,1]^d$, we have by \cite[Proposition 2]{Owe05} 
\begin{equation*}
\Delta (f;[0,1]^d)= \sum_{\y\in \mathcal{Y}} \Delta(f;[\y,\y_+]).
\end{equation*}
Define the variation over $\mathcal{Y}$ by 
\begin{equation*}
V_{\mathcal{Y}} (f;[0,1]^d)= \sum_{\y\in \mathcal{Y}} |\Delta(f;[\y,\y_+])|.
\end{equation*}
Then the \emph{Vitali variation} of $f$ over $[0,1]^d$ is defined by
\begin{equation*}
V (f;[0,1]^d)= \sup_{\mathcal{Y}\in \mathbb{Y}} V_{\mathcal{Y}} (f;[0,1]^d).
\end{equation*}
For a subset $u\subseteq \{1,\dots, d\}$, let 
\begin{equation*}
\Delta_u (f;[\abf,\bbf])= \sum_{v\subseteq u} (-1)^{|v|}f(\abf^v:\bbf^{-v}).
\end{equation*}
Let $\mathbf{0} = (0, \dots, 0) \in [0,1]^d$ and $\mathbf{1}=(1, \dots, 1)\in [0,1]^d$. Given a ladder $\Y$, there is a corresponding ladder $\mathcal{Y}_u = \{ \y^u:\mathbf{1}^{-u} \mid \y\in \Y \} $ on the $|u|$-dimensional face of $[0,1]^d$ consisting of points of the form $\x^{u}:\mathbf{1}^{-u}$ (we interpret $\mathcal{Y}_\emptyset $ as $\{\mathbf{1}\}$). The operation of the successor is also defined on $\mathcal{Y}_u$, and  again we have,
\begin{equation*}
\Delta_u (f;[0,1]^d)= \sum_{\y\in \mathcal{Y}_u} \Delta_u(f;[\y,\y_+]). 
\end{equation*}
Furthermore, we define
\begin{equation*}
V_{\mathcal{Y}_u} (f;[0,1]^d)= \sum_{\y\in \mathcal{Y}_u} |\Delta(f;[\y,\y_+])|,
\end{equation*}
which is the variation over the ladder $\Y_u$ of the restriction of $f$ to the face of $[0,1]^d$ specified by $u$.
The \emph{Hardy--Krause} variation is defined as 
\begin{equation*}
{HK}(f;[0,1]^d) = \sum_{\emptyset \neq u\subseteq \{1,\dots, d\}} \sup_{\mathcal{Y}\in \mathbb{Y}} V_{\mathcal{Y}_u} (f;[0,1]^d).
\end{equation*}
$\mathcal{HK}$ denotes the class of functions with bounded Hardy--Krause variation.
In words, the Hardy--Krause variation is the sum of the Vitali variations of all the restrictions of $f$ to those faces of $[0,1]^d$ adjecent to $\mathbf{1}$. \\

{\bf Leonov's result. } We follow \cite{AisDic15} and call a function  $f: [0,1]^d \to \R$ \emph{completely monotone} if the restriction $f_{\mid R}$ of $f$ to any axis-parallel box $R=[\abf,\bbf] \subseteq [0,1]^d$ of dimension $1 \leq s \leq d$ with $\abf\leq \bbf$ satisfies $\Delta^{(s)}(f_{\mid R}, [\abf,\bbf])\geq 0$.  The  $s$ in $\Delta^{(s)}$ is the dimension of $[\abf,\bbf]$ and marks that $f_{\mid R}$ is considered as a function of $s$ variables  when computing $\Delta^{(s)}$.  
We shall need the following result by Leonov \cite{leonov}: 
\begin{lem}[Leonov \cite{leonov}]\label{leo}
Any function of bounded Hardy--Krause variation can be written as the difference of two completely monotone functions.
\end{lem}

\subsection{$\mathcal{D}$-variation}
In the following, we recall the notion of variation introduced in \cite{PauSva14}.
Let $\mathcal{D}$ denote an arbitrary family of measurable subsets of $[0,1]^d$ with $\emptyset, [0,1]^d \in \mathcal{D}$. Let $\mathcal{S}(\mathcal{D})$ denote the corresponding vector space of simple functions 
\begin{equation*}
f=\sum_{i=1}^m \alpha_i \mathds{1}_{A_i}
\end{equation*}
where $\alpha_i \in \R$, $A_i\in \mathcal{D}$, and $m\in \N$. Note that the representation of $f$ is of course not unique.
We say that a set $A\subseteq [0,1]^d$ is an \emph{algebraic sum} of sets in $\mathcal{D}$ if there exist $A_1,\dots,A_m \in \mathcal{D} $ such that 
\begin{equation*}
\mathds{1}_A = \sum_{i=1}^n \mathds{1}_{A_i} - \sum_{i=n+1}^m \mathds{1}_{A_i}, 
\end{equation*}
and we define $\mathcal{A}$ to be the collection of algebraic sums of sets in $\mathcal{D}$.\\

Inspired by \cite{Har10}, we define the \emph{Harman complexity} ${h}(A) $ of a set $A\in \mathcal{A}$ with $A\neq [0,1]^d$ and $A\neq \emptyset$, as the minimal number $m$ such that there exists $A_1,\dots,A_m$ with 
\begin{align*}
\mathds{1}_A = \sum_{i=1}^n \mathds{1}_{A_i} - \sum_{i=n+1}^m\mathds{1}_{A_i}
\end{align*}
for some $n \in \{0, \dots, m\}$ and either $A_i \in \mathcal{D}$ or $[0,1]^d\backslash A_i\in \mathcal{D}$. Moreover, we define $h([0,1]^d)=h(\emptyset)=0$.\\

The definition of variation is given in two steps. First,  for $f\in \mathcal{S}(\mathcal{D})$, we define 
\begin{equation*}
\Vs {\mathcal{D}} (f) := \inf \bigg\{ \sum_{i=1}^m |\alpha_i| h_{\mathcal{D}}(A_i) \quad \bigg| \quad f=\sum_{i=1}^m \alpha_i \mathds{1}_{A_i}, \quad \alpha_i\in \R, A_i\in \mathcal{D}\bigg\}.
\end{equation*}
Second, let $\Vi {\mathcal{D}}$ be the collection of all measurable functions $f :[0,1]^d \to \R$ for which there exists a sequence of $f_i \in \mathcal{S}(\mathcal{D})$ that converges  to $f$ in the supremum norm $|\cdot|_\infty$.

\begin{defi}[\cite{PauSva14}, Definition 3.2] \label{variation}
We define the $\mathcal{D}$-variation of $f\in \Vi {\mathcal{D}}$ as
\begin{equation*}
\Va {\mathcal{D}}(f)= \inf \, \Big\{ \, \liminf_i \Vs {\mathcal{D}} (f_i) \quad \Big| \quad f_i \in \mathcal{S}(\mathcal{D}), \, \lim_i |f-f_i|_\infty =0 \, \Big\}
\end{equation*}
and set $\Va {\mathcal{D}}(f)=\infty $ if $f\notin \Vi {\mathcal{D}}$. The space of functions of bounded $\mathcal{D}$-variation is denoted by
\begin{equation*}
\Vf {\mathcal{D}} = \Big\{f\in \Vi {\mathcal{D}} \quad \Big| \quad \Va {\mathcal{D}}(f)< \infty \Big\}.
\end{equation*}
\end{defi}
Among the classes of sets $\mathcal{D}$ which are of particular interest are the class $\mathcal{K}$ of convex sets and the class $\mathcal{R}^\ast$ of axis parallel boxes containing $\mathbf{0}$ as a vertex.
In the following we recall the most important properties of this notion.
\begin{prop} \label{properties}\textcolor{white}{lala}\\
(i) $\Vi {\mathcal{D}}$ and $\Vf {\mathcal{D}}$ are vector spaces.
In particular, $\Va {\mathcal{D}}$ defines a semi-norm on $\Vf {\mathcal{D}}$. \\
(ii) $\Vi {\mathcal{D}}$ is closed under limits in the supremum-norm. We have the following lower semi-continuity: if $|f-f_i|_\infty \to 0$ then
$\Va {\mathcal{D}}(f) \leq \liminf_i \Va {\mathcal{D}}(f_i)$.\\
(iii) If $\mathcal{D}$ is closed under intersection, then
$\Vf {\mathcal{D}}$ is closed under multiplication, and
\begin{align}\label{multibound}
\Va {\mathcal{D}}(fg) \leq 3 \Va {\mathcal{D}}(f) \Va {\mathcal{D}}(g) + \inf|f| \Va {\mathcal{D}}(g) + \inf|g|\Va {\mathcal{D}}(f).
\end{align}
\end{prop}

\begin{Proof}
(i) is \cite[Proposition 3.5]{PauSva14}, (ii) is \cite[Proposition 3.6]{PauSva14}, (iii) is \cite[Theorem 3.7]{PauSva14}.
\end{Proof}


\section{Borel measurability of functions of bounded Hardy--Krause variation}
\label{sec3}

The aim of this section is to give an independent proof that every function of bounded Hardy--Krause variation is Borel measurable. This fact also plays a key role in the equivalence of Hardy--Krause variation and $\mathcal{R}^*$-variation, which will be stated and proved in the subsequent section.

\begin{thm} \label{thm:measure}
Every function of bounded Hardy--Krause variation is Borel measurable. More precisely, every real-valued function on $[0,1]^d$ which has bounded HK-variation is \\
$\left([0,1]^d,\mathcal{B}\left([0,1]^d\right)\right)-\left(\mathbb{R},\mathcal{B}(\mathbb{R})\right)$-measurable.
\end{thm}

We have looked in the literature very carefully, but have not found anywhere the fact that finite HK-variation implies Borel measurability. It is remarkable that such a fundamental property of functions of bounded HK-variation has not been investigated before. However, the proof is far from being trivial (see below).  Recall from Lemma \ref{leo} that a function of bounded HK-variation decomposes into a difference of two completely monotone functions (see also \cite{AisDic15}), so the assertion of Theorem \ref{thm:measure} follows from a similar result for completely monotone functions, stated in Theorem \ref{thm:complete} below. The fact that completely monotone functions are Borel measurable also seems to be new. Note that coordinatewise monotonicity is \emph{not} sufficient for a multivariate function to be Borel measurable. For a two-dimensional counterexample, define $f(x,y)=0$ for $x+y<1$, $f(x,y)=1$ for $x+y>1$, and for $x+y=1$ set $f(x,y)=1/2$ for $x \in E$ and $f(x,y)=0$ 
otherwise, where $E \subset [0,1]$ is not Borel measurable. Then $f^{-1} (\{1/2\})$ is not in $\mathcal{B}\left([0,1]^2\right)$.\\

On the other hand, coordinatewise monotonicity \emph{is} actually sufficient for \emph{Lebesgue} measurability of a multivariate function. That means, a function which is coordinatewise monotone is $\left([0,1]^d,\mathcal{L}\left([0,1]^d\right)\right)-\left(\mathbb{R},\mathcal{B}(\mathbb{R})\right)$-measurable, where $\mathcal{L}$ is the Lebesgue sigma-field (the completion of the Borel sigma-field). A possible proof goes as follows. We use induction on $d$. The case $d=1$ is trivial. Now let $f$ be a function of $d$ variables which is increasing in every coordinate. For fixed $a \in \mathbb{R}$, define $g(x_1, \dots, x_{d-1}) = \sup \{y \in [0,1]:~f(x_1, \dots, x_{d-1},y) \leq a\}$, where the supremum of the empty set is understood to be zero. Then $g$ is monotonic decreasing, and, by the induction hypothesis, Lebesgue measurable. Thus the set $A = \{(x_1, \dots, x_{d-1},y):~g(x_1, \dots, x_{d-1}) < y\}$ is also Lebesgue measurable. Moreover, the set $B=\{(x_1, \dots,x_{d-1},y):~g(x_1, \dots, x_{d-1})=y\}$ is of measure zero (by Fubini's theorem). Now the set $\{f > a\}$ differs from $A$ only by a subset of $B$, which has measure zero, and hence is Lebesgue measurable.\\

The remarks above show that these measurability issues are rather delicate, and should be treated very carefully (as a deterrent example cf. \cite{idcz}, where it is proved that every multivariate, coordinatewise monotonic function is ``measurable'', without any mention in the whole paper which kind of measurability is actually meant -- in fact the author talks about Lebesgue measurability, but careless readers may easily be misled).

\begin{thm} \label{thm:complete}
Every completely monotone and real-valued function on $[0,1]^d$ is\\ $\left([0,1]^d,\mathcal{B}\left([0,1]^d\right)\right)-\left(\mathbb{R},\mathcal{B}(\mathbb{R})\right)$-measurable.
\end{thm}

\begin{Proof}[Proof of Theorems \ref{thm:measure} and \ref{thm:complete}.]
By Lemma \ref{leo}, a function of bounded HK-variation can be written as the difference of two completely monotone functions. Thus Theorem \ref{thm:measure} is a consequence of Theorem \ref{thm:complete}, and in the sequel we will assume that $f$ is a completely monotone function.\\

We proceed by induction on the number of variables $d$. In the case $d=1$ complete monotonicity reduces to (ordinary) monotonicity, and the Borel measurability of monotonic functions in one variable is a classical result. This proves the initial step of the induction.\\

Now we assume that the induction hypothesis holds for all completely monotone functions which have less than $d$ variables, and assume that $f$ is a completely monotone function on $[0,1]^d$. It is a well-known fact that all the discontinuities of a completely monotone function lie on an at most countable set of hyperplanes of dimensions $d-1$, all of which are parallel to the coordinate axes (this fact was probably first noted by Young and Young \cite{young}, and rediscovered by Antosik \cite{antosik}). As a consequence, roughly speaking, $f$ decomposes into a continuous part (which is measurable by continuity) and into countably many lower-dimensional functions (which are measurable by the induction hypothesis). However, this argument has to be carried out 
very carefully; all the details are given below.\\

We write $H_1, H_2, \dots$ for the collection of $(d-1)$-dimensional hyperplanes where the discontinuities of $f$ are situated, and we set
\begin{equation} \label{ddef}
D = [0,1]^d \backslash \left( \bigcup_{k=1}^\infty H_k \right).
\end{equation}
Since $f$ is assumed to be completely monotone, by definition it is monotonically increasing. Thus there exists a number $m$ such that $m < f(\mathbf{0}) \leq f(\mathbf{x})$ for all $\mathbf{x} \in [0,1]^d$. For $k \geq 1$, we define
$$
f_k (\mathbf{x}) = \left\{\begin{array}{ll} f(\mathbf{x}) & \textrm{for $\mathbf{x} \in H_k$},\\ m & \textrm{otherwise}. \end{array}\right.
$$
Let $k$ be given. Then there exist an index $i \in \{1, \dots, d\}$ and a number $a \in [0,1]$ such that the hyperplane $H_k$ consists of all the points $\left\{\mathbf{x} = \left(x^{(1)},\dots,x^{(d)} \right) \in [0,1]^d:~x^{(i)} = a \right\}$. Furthermore, the $d$-variate function $f_k$ induces in a natural way a $d-1$-variate function $\hat{f}_k$ by the relation
\begin{eqnarray}
& & \hat{f}_k \left(x^{(1)}, ~x^{(2)}, ~\dots, x^{(i-1)}, ~~~~~~ x^{(i+1)},~ \dots,~ x^{(d-1)}, ~x^{(d)}\right) \nonumber\\
& = & f_k \left(x^{(1)}, ~x^{(2)}, ~\dots, ~x^{(i-1)}, ~a, ~x^{(i+1)}, ~\dots, ~x^{(d-1)}, ~x^{(d)} \right). \label{hatf}
\end{eqnarray}
By the definition of complete monotonicity, the function $\hat{f}_k$ is a $(d-1)$-variate completely monotone function. Consequently, by the induction hypothesis, $\hat{f}_k$ is Borel measurable on $[0,1]^{d-1}$. Thus the preimage of a Borel set of $\mathbb{R}$ under $f_k$ consists of
\begin{itemize}
 \item The part contained in $H_k$, which is the cross-product of a Borel set of $[0,1]^{d-1}$ and of a one-point set (the point $a$ in equation \eqref{hatf}), and which consequently is $\mathcal{B}\left([0,1]^d\right)$-measurable.
 \item Possibly additionally the whole set $[0,1]^d \backslash H_k$, which is also measurable.
\end{itemize}
Thus for every $k$ the function $f_k$ is a measurable function from $\left([0,1]^d,\mathcal{B}\left([0,1]^d\right)\right)$ to $\left(\mathbb{R},\mathcal{B}(\mathbb{R})\right)$.\\

Next we define a function $g$ by setting
$$
g (\mathbf{x}) = \left\{\begin{array}{ll} f(\mathbf{x}) & \textrm{for $\mathbf{x} \in D$},\\ m & \textrm{otherwise}, \end{array}\right.
$$
where $D$ is the set from \eqref{ddef}. We want to show that for every given $b \in \mathbb{R}$ the preimage of $(-\infty,b)$ under $g$ is measurable. Then, since the Borel sigma-field on $\mathbb{R}$ is generated by the collection of sets $\{(-\infty,b), ~b \in \mathbb{R}\}$, the function $g$ is a measurable function from $\left([0,1]^d,\mathcal{B}\left([0,1]^d\right)\right)$ to $\left(\mathbb{R},\mathcal{B}(\mathbb{R})\right)$ (see \cite[Theorem 1.41]{yeh}). Thus let $b \in \mathbb{R}$ be fixed, and set $B = (-\infty,b)$. If $b \leq m$, then by construction the set $g^{-1}(B)$ is the empty set (which is measurable). If $b > m$, then by construction we have
$$
\left([0,1]^d \backslash D\right) \subset g^{-1} (B).
$$
Furthermore, since $f$ and $g$ coincide on $D$, we have
\begin{equation} \label{hatfb}
g^{-1}(B) = f^{-1}(B) \cup \left([0,1]^d \backslash D\right).
\end{equation}
Now assume that $\mathbf{x} \in D$, and that $\mathbf{x} \in f^{-1}(B)$. Then there exists a number $y \in B$ such that $y = f(\mathbf{x})$. By construction the function $f$ is continuous in $\mathbf{x}$. Note that the set $B$ is open, which implies that there exists a $\delta>0$ such that a $\delta$-neighborhood around $y$ is also contained in $B$. Accordingly, by the definition of continuity, there exists an $\varepsilon>0$ such that all elements of $[0,1]^d$ which are contained in an $\varepsilon$-neighborhood of $\mathbf{x}$ are mapped by $f$ into the $\delta$-neighborhood of $y$. Thus there exists an open set $N_{\mathbf{x}} \subset [0,1]^d$ containing $\mathbf{x}$ such that $f(N_{\mathbf{x}}) \subset B$. It is easily verified that we have
$$
\bigcup_{\mathbf{x} \in D,~f(\mathbf{x})  \in B} N_{\mathbf{x}} \subset f^{-1} (B) \subset \left(\bigcup_{\mathbf{x} \in D,~f(\mathbf{x})  \in B} N_{\mathbf{x}}\right) \cup \left([0,1]^d \backslash D\right).
$$
Thus by \eqref{hatfb} we have
\begin{equation} \label{hatfb2}
g^{-1}(B) = \left(\bigcup_{\mathbf{x} \in D,~f(\mathbf{x})  \in B} N_{\mathbf{x}}\right) \cup \left([0,1]^d \backslash D\right).
\end{equation}
The set on the right-hand side of \eqref{hatfb2} is the union of 
\begin{itemize}
 \item a union of open sets (which itself is also open, and consequently Borel measurable), and of
 \item a countable union of hyperplanes (which also is Borel measurable).
\end{itemize}
Thus $g^{-1}(B) \in \mathcal{B}\left([0,1]^d\right)$, which proves that $g$ is a measurable function from $\left([0,1]^d,\mathcal{B}\left([0,1]^d\right)\right)$ to $\left(\mathbb{R},\mathcal{B}(\mathbb{R})\right)$.\\

Thus we have established that all the functions $g$ and $f_k,~k \geq 1$, are measurable. By construction we have
$$
f(\mathbf{x}) = \sup_{N \geq 1} \max \Big\{g(\mathbf{x}),f_1(\mathbf{x}),f_2(\mathbf{x}), \dots, f_N(\mathbf{x})\Big\}, \qquad \mathbf{x} \in [0,1]^d.
$$
The supremum of measurable functions is itself measurable (see \cite[Theorem 4.22]{yeh}). Thus we have established that $f$ is a measurable function from $\left([0,1]^d,\mathcal{B}\left([0,1]^d\right)\right)$ to $\left(\mathbb{R},\mathcal{B}(\mathbb{R})\right)$, which proves the theorem.
\end{Proof}


\section{Equivalence of Hardy--Krause and $\mathcal{R}^*$-variation}
\label{sec4}

The aim of this section is to show that $\mathcal{D}$-variation with  $\mathcal{D}=\mathcal{R}^*$ coincides with Hardy--Krause variation.
The following was already shown in \cite{PauSva14}. 
\begin{thm}[\cite{PauSva14}]\label{Oneway}
$\mathcal{HK} \cap \Vi {\mathcal{R}^*}=\Vf {\mathcal{R}^*}$ and  $HK(f;[0,1]^d) = \Va{\mathcal{R}^*}(f)$ whenever $f\in \Vf {\mathcal{R}^*}$.
\end{thm} 
We shall show the following theorem. 
\begin{thm}\label{equal}
Every function of bounded Hardy--Krause variation can be uniformly approximated by a sequence of simple functions from $\mathcal{S}(\mathcal{R}^*)$, i.e.\  $\mathcal{HK}\subseteq \Vi {\mathcal{R}^*}$. 
\end{thm}
Combining this with Theorem \ref{Oneway} yields:
\begin{cor}\label{variationequality}
We have $\mathcal{HK} = \Vf {\mathcal{R}^*}$, and for any $f:[0,1]^d \to \R$ we have
\begin{equation*}
HK(f;[0,1]^d) = \Va{\mathcal{R}^*}(f).
\end{equation*}
\end{cor}
Since the limit of a sequence of measurable functions is again measurable, we immediately obtain Theorem \ref{thm:measure} as a corollary.
\begin{cor}\label{measureable}
Every function of bounded Hardy--Krause variation is Borel measureable.
\end{cor}

\begin{rem}
Corollary \ref{variationequality} shows that ${\mathcal{R}^*}$-variation yields an alternative way of constructing Hardy--Krause variation. This is convenient in some situations, as illustrated by Corollary \ref{measureable}: Measurability is obvious from the definition of ${\mathcal{R}^*}$-variation, whereas the proof based on the classical definition is involved and relies on results of earlier papers on the points of discontinuity of a function of bounded variation. Other properties, such as the fact that if $f> \delta >0 $ has bounded Hardy--Krause variation, then so does $1/f$, are easily shown with the classical definition, but it is not clear to the authors how to obtain this fact directly from the definition of ${\mathcal{R}^*}$-variation.
\end{rem}
\begin{rem}
Corollary \ref{variationequality} also shows that $\mathcal{D}$-variation is a quite general concept in the sense that the spaces $\Vf {\mathcal{D}}$ of functions of bounded $\mathcal{D}$-variation are rather large and contain many interesting functions. Indeed, if $\mathcal{R}^* \subseteq \mathcal{D}$ then $\mathcal{HK} \subseteq \Vf{\mathcal{D}}$. In particular, the space of functions of bounded $\mathcal{K}$-variation contains $\mathcal{HK}$, but is known to be strictly larger. This was not at all clear to the authors in \cite{PauSva14}.
\end{rem}

\begin{rem}
A Koksma-Halwka inequality for general measures on $[0,1]^d$, which is the main result in G\"{o}tz \cite{gotz} and is also stated as Theorem 1 in \cite{AisDic15}, follows directly from Corollary \ref{variationequality} and the Koksma-Hlawka inequality \cite[Thm. 4.3]{PauSva14}. 
\end{rem}

Before we prove Theorem \ref{equal}, we slightly extend the notation introduced in Section \ref{sec2} and present an important observation on completely monotone functions. 
If $v \subseteq \{1,\ldots, d\}$ and $\mathbf{a}, \mathbf{b} \in [0,1]^d$, then $[\mathbf{a},\mathbf{b}]^v$ denotes the box 
\begin{equation*}
[\mathbf{a},\mathbf{b}]^v=\{ \mathbf{x} \in \R^d \mid \forall i\in v : a_i\leq x_i \leq b_i, \forall i \notin v : a_i\leq x_i < b_i\}.
\end{equation*}
If $a_i = b_i $ for some $i\notin v$, then $[\mathbf{a},\mathbf{b}]^v$ should be interpreted as the empty set.
The set of axis-parallel rectangles containing $\mathbf{0}$ is then given by
\begin{equation*}
\mathcal{R}^* =\{[\mathbf{0},\mathbf{a}]^v \mid \mathbf{a}\in [0,1]^d, v\subseteq \{1,\ldots,d\}\}.
\end{equation*}
Given a ladder $\mathcal{Y}$, we define a partial ordering of the pairs $(\mathbf{y},v)$ with $\mathbf{y}\in \bigcup_{u \subseteq  \{1,\dots,d\}}\mathcal{Y}_u$ and $v\subseteq \{1,\dots,d\}$ by declaring $(\mathbf{y},v)\leq  (\mathbf{z},w)$ if $\mathbf{y}\leq \mathbf{z}$ and $\mathbf{y}\neq \mathbf{z}$ or if $\mathbf{y}=\mathbf{z}$ and $v \subseteq w$. 
We denote by $C(\mathbf{y},v)$ the face 
\begin{equation*}
[\mathbf{0},\mathbf{y}]^v\backslash \bigcup_{(\mathbf{z},w)<(\mathbf{y},v)} [\mathbf{0},\mathbf{z}]^{w}.
\end{equation*}
Intuitively, $C(\mathbf{y},v)$ denotes the $d-| v|$-dimensional face of the subrectangle formed by the ladder whose maximal vertex is $\mathbf{y}$ and whose spanning edges specified by $v$ have maximal vertex $\mathbf{y}$; see Figure \ref{cyv} (left).

\begin{center}
\begin{figure}[h!]
\centering
\begin{tikzpicture}[scale=0.5]

\draw[gray, very thin] (0,0) -- (0,4) -- (4,4) -- (4,0) -- (0,0);
\draw[gray, very thin] (2,0) --(2,4);
\draw[gray, very thin] (0.7,0) --(0.7,4);
\draw[gray, very thin] (3.6,0) --(3.6,4);
\draw[gray, very thin] (0,2) --(4,2);
\draw[gray, very thin] (0,1) --(4,1);
\draw[gray, very thin] (0,2.8) --(4,2.8);

\draw[thick] (2,2) -- (2,1);

\node at (0.7,-0.4) {{\scriptsize $y_1^1$}};
\node at (2,-0.4) {{\scriptsize $y_2^1$}};
\node at (3.6,-0.4) {{\scriptsize $y_3^1$}};

\node at (-0.4,1) {{\scriptsize $y_1^2$}};
\node at (-0.4,2) {{\scriptsize $y_2^2$}};
\node at (-0.4,2.8) {{\scriptsize $y_3^2$}};


\draw[dashed] (8,0) -- (11,0) -- (11,3) -- (8,3) -- (8,0);
\draw[dashed] (8,0) -- (9,1) -- (12,1) -- (12,4) -- (9,4) -- (9,1);
\draw[dashed] (8,3) -- (9,4);
\draw[dashed] (11,3) -- (12,4);
\draw[dashed] (11,0) -- (12,1);
\draw[thick] (9,4) -- (12,4);

\node at (8,-0.4) {{\scriptsize $(0,0,0)$}};
\node at (12.4 , 4.4) {{\scriptsize $(1,1,1)$}};

\end{tikzpicture}
\caption{Left: A grid ladder (gray) and the face $C( (y_2^1, y_2^2) ,\{1\})$ (bold). Right: The cube $[0,1]^3$ (dashed) and the face $F_1$ (bold) -- this will be defined and used at the end of this section.} \label{cyv}
\end{figure}
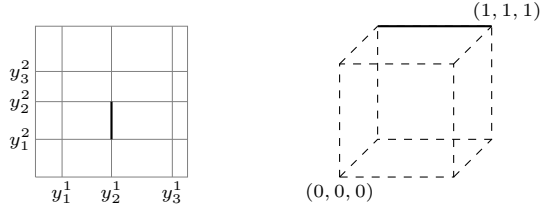
\end{center}
\vspace{-15pt}

Moreover, for $\mathbf{x},\mathbf{y},\mathbf{a} \in [0,1]^d$ and $i,j \in \{1,\ldots, d\}$ with $i\neq j$ we write $\mathbf{y}^i:\mathbf{x}^j:\mathbf{a}^{-i\cup j}$ for the point with $i$-th coordinate $y_i$,  $j$-th coordinate  $x_j$, and all other coordinates equal to those of $\mathbf{a}$.

\begin{lem}\label{help}
Let $f$ be a completely monotone function and let $\mathbf{x},\mathbf{a} \in [0,1]^d$ and  $ i,j\in \{1,\dots , d\}$, $i\neq j$. Then
\begin{align*}
|f( \mathbf{x}^j: \mathbf{a}^{-j}) -  f(\mathbf{a}) | {}&\leq  |f(\mathbf{1}^{i}: \mathbf{x}^j : \mathbf{a}^{-i\cup j})  - f(\mathbf{1}^{i} : \mathbf{a}^{-i })|\\
&\leq  |f(\mathbf{x}^j : \mathbf{1}^{- j})  - f( \mathbf{a}^j : \mathbf{1}^{-j})|.
\end{align*}
\end{lem}

\begin{Proof}
We use the complete monotonicity of $f$ restricted to the $2$-dimensional plane $\{\mathbf{z}\in [0,1]^d: \mathbf{z}^{-i\cup j}= \mathbf{a}^{-i\cup j}\}$.  If $a_j \leq x_j$, we get
\begin{equation*}
0 \leq f(\mathbf{a}) + f(\mathbf{1}^{i}: \mathbf{x}^j : \mathbf{a}^{-i\cup j}) - f( \mathbf{x}^j: \mathbf{a}^{-j}) - f(\mathbf{1}^{i}: \mathbf{a}^{-i}),
\end{equation*}
which, together with the monotonicity on $1$-dimensional spaces, yields 
\begin{equation*}
0\leq f( \mathbf{x}^j: \mathbf{a}^{- j}) -  f(\mathbf{a})  \leq  f(\mathbf{1}^{i}: \mathbf{x}^j : \mathbf{a}^{-i\cup j})  - f(\mathbf{1}^{i} :  \mathbf{a}^{-i}).
\end{equation*}
If $a_j \geq x_j$, all inequalities are reversed, so in both cases, we can deduce
\begin{equation*}
|f( \mathbf{x}^j: \mathbf{a}^{-j}) -  f(\mathbf{a}) | \leq  |f(\mathbf{1}^{i}: \mathbf{x}^j : \mathbf{a}^{-i\cup j})  - f(\mathbf{1}^{i} : \mathbf{a}^{-i})|.
\end{equation*}
The second inequality follows by repeated use of the first one.
\end{Proof}

\begin{Proof}[Proof of Theorem \ref{equal}.]
By Lemma \ref{leo}, we may assume that $f$ is completely monotone. 
We prove the theorem in two steps. We first study the case $d=1$, and then we apply this result to approximate functions on $[0,1]^d$ for general $d$ as well.\\

First let $d=1$. A completely monotone function $f:[0,1]\to \R$ is a bounded increasing function. In particular, we may choose a simple function $g_n$ with $|f-g_n|_\infty \leq 1/n$.
One can construct $g_n$ as follows: Choose a partion $\inf f=t_1 < \dotsm < t_{N+1}=\sup f$ of the interval $[\inf f,\sup f]$ such that $|t_l -t_{l+1}| \leq 1/n$ for all $l$. Then choose $y_l\in f^{-1}(t_l)$ if such a point exists. Otherwise $ f^{-1}(t_l)=\emptyset$, which means that $f$ has a jump from values smaller than $t_l$ to values larger than $t_l$ at some point. We denote this point by $y_l$. We may take $y_1=0$, $y_{N+1}=1$. (Some of the $y_l$ may be equal. In this case, we throw away multiple points so that $y_l < y_{l+1}$ for all $l$.)\\

On each interval $(y_l,y_{l+1})$, $f$ is increasing and its function values lie in $[t_l , t_{l+1}]$. Therefore, if we choose  a point $z_l \in (y_l,y_{l+1})$, then for any $ x \in (y_l,y_{l+1})$ we have $|f(x)-f(z_{l})| \leq |t_{l+1}-t_l| \leq 1/n$. 
Thus, if we define $g_n$ such that $g_n(x)=f(z_l)$ for $x\in (y_l, y_{l+1})$ and $g_n(x)=f(y_l)$ for $x=y_l$, then $|g_n - f|_\infty \leq 1/n$.\\

We need to show that $g_n \in \mathcal{S}(\mathcal{R}^*)$. To achieve this, we define a new step function $f_n$ as a sum of contributions from the half open intervals $[0, y_l)$ and the closed intervals $[0,y_l]$. The contribution to $f_n$ from the half open intervals are $f(z_{l-1}) - f(y_{l})$, and  $f(y_l) - f(z_l)$ from the closed intervals (with the exception of the interval $[0,1]$ which contributes $f(1)$). In this way, $f_n(x)$ is written as a telescoping sum and  evaluates to either $f(y_l)$ or $f(z_l)$ depending on whether $x=y_l$ or $x \in (y_l, y_{l+1})$. Thus $f_n(x)$ agrees with $g_n(x)$ for all $x \in [0,1]$. 
Using the heavy notation introduced above, we can write $f_n$ as
\begin{equation}\label{1dapprox}
f_n =  \sum_{l=1}^{N+1} \sum_{v\subseteq \{1\}} \alpha_{l,v} \mathds{1}_{[0,y_l]^{v}}.
\end{equation}
Here $\alpha_{N+1,\{1\}}=f(1)$, and for $(l,v)< (N+1,\{1\})$, we define $\alpha_{l,v}$ as follows: Let $z_{l,\{1\}}=y_l$ and $z_{l,\emptyset}=z_{l-1}$. 
Let $\mathcal{Y}$ be the ladder formed by $y_1,\dots , y_N$ and $\tilde{\mathcal{Y}}$ the ladder formed by the $z_{l,v}$ with $(l,v)< (N+1,\{1\})$. Then
\begin{equation*}
\alpha_{l,v} = -{\Delta}(f;z_{l,v},(\tilde{z}_{l,v})_+ ) = f(z_{l,v})-f((\tilde{z}_{l,v})_+).
\end{equation*}
Here $\tilde{y}_+$ indicates that the successor of $y$ is computed with respect to the ladder $\tilde{\mathcal{Y}}$.  With this definition, $f_n(x)$ agrees with $g_n(x)$ for all $x \in [0,1]$, which concludes the argument in the one-dimensional case. \\

This idea can be extended to the multi-dimensional case by replacing the half open intervals by partially open axis aligned boxes to which we attach alternating sums of function values. 
To see this, let $d>1$ and consider a completely monotone function $f$ on $[0,1]^d$. Let $f^i$, $i=1,\dots,d$, denote the restriction of $f$ to the 1-dimensional face  $F_i =\{\mathbf{x}\in [0,1]^d \mid \mathbf{x}^{-i}=\mathbf{1}^{-i}\}$ of $[0,1]^d$; see Figure \ref{cyv} (right) for an illustration of $F_i$. We apply the result for $d=1$ to each $f^i$ and
 define ladders $\mathcal{Y}^i$ and $\tilde{\mathcal{Y}}^i$ consisting of points $\mathbf{y}_l^i$ and $\mathbf{z}_{l,v}^i$, respectively, such that we have an approximation $f_n^i$ of the form \eqref{1dapprox}, i.e.
\begin{equation*}
f_n^i =  \sum_{l=1}^{N_i+1} \sum_{v\subseteq \{1\}}\alpha_{l,v}^i \mathds{1}_{[\mathbf{0},\mathbf{y}_l^i]^{v}}
\end{equation*}
and $|f^i- f^i_n|_\infty \leq 1/n $.\\

Form the ladders $\mathcal{Y}= \prod_i \mathcal{Y}^i$ and $\tilde{\mathcal{Y}}= \prod_i \tilde{\mathcal{Y}}^i$. For $\mathbf{y}=(y^1_{l_1},\dots,y^d_{l_d}) \in \mathcal{Y}_u$, let $\mathbf{z}_{ y,v}$ be the point whose $i$-th coordinate is $z_{l_i, \{1\}}^i=y^i_{l_i}$ if $i\in v$, and $z_{l_i, \emptyset}^i$ otherwise. Then $\mathbf{z}_{y,v}\in C(\mathbf{y},v)$. 
Consider the simple function
\begin{equation*}
f_n(\mathbf{x}) = \sum_{u\subseteq \{1,\dots,d\}} ~\sum_{\mathbf{y}\in \mathcal{Y}_u}~ \sum_{v\subseteq \{1,\dots,d\}} \alpha_{u,\mathbf{y},v} \mathds{1}_{[\mathbf{0},\mathbf{y}]^{v}}(\mathbf{x})
\end{equation*}
where 
\begin{equation*}
\alpha_{u,\mathbf{y},v}=(-1)^{|u\cup(-v)|}\Delta_{u\cup (-v)}(f,\mathbf{z}_{y,v},\tilde{\mathbf{z}}_{y,v}^+).
\end{equation*}
Here $\tilde{\mathbf{y}}^+$ is the successor of $\mathbf{y}$ in its ladder. \\

Theorem \ref{equal} is proved if we can show that $|f_n - f|_\infty \leq d/n$.  
Note that $f_n$ is constantly equal to $f(\mathbf{z}_{y_0,v_0})$ on each $ C(\mathbf{y}_0,v_0) $ because  for $\mathbf{x}\in C(\mathbf{y}_0,v_0) $ we have
\begin{align*}
f_n(\mathbf{x}) {}&= \sum_{u\subseteq \{1,\dots,d\}} ~\sum_{\mathbf{y}\in \mathcal{Y}_u}~ \sum_{v\subseteq \{1,\dots,d\}} \alpha_{u,\mathbf{y},v} \mathds{1}_{[\mathbf{0},\mathbf{y}]^{v}}(\mathbf{x})\\
&=\sum_{u\subseteq \{1,\dots,d\}} ~\sum_{\mathbf{y}\in {\mathcal{Y}}_u}~ \sum_{v\subseteq \{1,\dots,d\}}  (-1)^{|u\cup(-v)|}\Delta_{u\cup (-v)}(f,\mathbf{z}_{y,v} ,\mathbf{\tilde{z}}_{y,v}^+)\mathds{1}_{[\mathbf{0},\mathbf{z}_{y,v}]}(\mathbf{z}_{{y}_0,v_0})\\
&=\sum_{u\subseteq \{1,\dots,d\}} ~\sum_{\mathbf{z}\in \tilde{\mathcal{Y}}_u} (-1)^{|u|}\Delta_{u}(f,\mathbf{z} ,\mathbf{\tilde{z}}_+)\mathds{1}_{[\mathbf{0},\mathbf{z}]}(\mathbf{z}_{{y}_0,v_0})\\
&=f(\mathbf{z}_{{y}_0,v_0}),
\end{align*}
where the last equality uses \cite[Proposition 6]{Owe05}.
Observe also that $f_n=f_n^i$ on each of the faces $F_i$.  Let $\mathbf{x}\in C({\mathbf{y},v})$ be given and write $\mathbf{q}=\mathbf{z}_{ y,v}$. Then 
\begin{align*}
|f(\mathbf{x})-f_n(\mathbf{x})|{}& = |f(\mathbf{x})-f(\mathbf{q})| \\
{}&\leq \sum_{i=1}^d  |f(x_1,\dots,x_i,q_{i+1},\dots, q_d) -  f(x_1,\dots , x_{i-1},q_{i},\dots , q_d)| \\
&\leq \sum_{i=1}^d |f(1,\dots,1,x_i,1\dots 1 ) -  f(1,\dots , 1,q_{i}, 1 \dots , 1)|\\
&= \sum_{i=1}^d |f^i(x_i ) -  f^i_n(q_{i} )|\\
& \leq d /n,
\end{align*}
where the second inequality follows from Lemma \ref{help} and the second equality used  that $ f_n^i(\mathbf{z}_{l,v}^i)=f(\mathbf{z}_{l,v}^i)$ by construction. 
\end{Proof}


\section{Algebraic structure of $\Vf {\mathcal{D}}$}
\label{sec5}

Finally, we consider the algebraic structure of the function space $\Vf {\mathcal{D}}$. We will assume throughout this section that $\mathcal{D} $ is closed under intersection, i.e. $D_1,D_2\in \mathcal{D}$ implies $D_1\cap D_2\in \mathcal{D}$.

\subsection{Algebraic Structure of $\Vf {\mathcal{D}}$}
Analogous to \cite{BluTic89} we define for $f \in \Vf {\mathcal{D}}$ and $\sigma > 0$
\begin{equation*}
\| f \| = \| f \|_{\infty} + \sigma \Va {\mathcal{D}}(f),
\end{equation*}
which is a norm on $\Vf {\mathcal{D}}$. In the following we show that for $\sigma \geq 3 $, ($\Vf {\mathcal{D}}$, $\| \cdot \|$) is a commutative Banach algebra with respect to pointwise multiplication.

\begin{lem} \label{complete}
The norm $\| \cdot \|$ is complete.
\end{lem}

\begin{Proof}
Let $(f_i)$ be a Cauchy sequence in $(\Vf {\mathcal{D}}, \| \cdot \|)$. Then $\|f_i-f_j\| < \varepsilon$ implies that
\begin{equation*}
\|f_i - f_j \|_{\infty} < \varepsilon \ \ \ \ \text{ and } \ \ \ \ \sigma \Va {\mathcal{D}}(f_i - f_j) < \varepsilon,
\end{equation*}
because both summands are nonnegative, since $\sigma>0, \| \cdot \|_{\infty}$ is a norm and $\Va {\mathcal{D}} $ is a seminorm. This implies that $(f_i)$ is a Cauchy sequence with respect to the supremum norm and hence converges uniformly to some $f \in \Vi {\mathcal{D}}$. 
Now choose a subsequence $i_k$ such that $\Va {\mathcal{D}} (f_{i_k}-f_{i_k+j}) < 1/2^k$ for all $j>0$. Then $(f_{i_k})_{k\geq 1}$ is also a Cauchy sequence with respect to the supremum norm and converges uniformly to $f$. The semi-continuity (see \cite[Proposition 3.6]{PauSva14}) yields that
\begin{equation*}
\Va {\mathcal{D}}(f-f_i) \leq \liminf_k \Va {\mathcal{D}}(f_{i_k}-f_i),
\end{equation*}
which is smaller than some given $\varepsilon >0$ if $i$ is sufficiently large. Thus, $f \in \Vf {\mathcal{D}}$.
\end{Proof}

The next thing is to show that the Banach space norm is submultiplicative.

\begin{lem} \label{submult}
Assume $\mathcal{D}$ is closed under intersections. Let $f,g \in (\Vf {\mathcal{D}}, \| \cdot \|)$ and let $\sigma \geq 3$. Then
$\| f g \| \leq \|f\| \|g\|$.
\end{lem}

\begin{Proof}
This follows via a direct calculation from \eqref{multibound}:
\begin{align*}
\| fg \| &= \|fg\|_{\infty} + \sigma \Va {\mathcal{D}} (fg) \\
		&\leq \|fg\|_{\infty} + \sigma \inf |f| \Va {\mathcal{D}}(g) +\sigma \inf |g| \Va {\mathcal{D}}(f)  + 3 \sigma \Va {\mathcal{D}}(f) \Va {\mathcal{D}}(g) \\
		& \leq \| f \|_{\infty} \| g \|_{\infty} +\sigma \|f\|_{\infty} \Va {\mathcal{D}}(g) + \sigma \|g\|_{\infty} \Va {\mathcal{D}}(f) + \sigma^2 \Va {\mathcal{D}}(f) \Vf {\mathcal{D}}(g) \\
	&= \|f \| \|g\|.
\end{align*}
\end{Proof}

\begin{thm}
If $\mathcal{D}$ is closed under intersections and $\sigma \geq 3$, then
($\Vf {\mathcal{D}}$, $\| \cdot \|$) is a commutative Banach algebra with respect to pointwise multiplication.
\end{thm}

\begin{Proof}
By Proposition \ref{properties} and Lemma \ref{complete}, ($\Vf {\mathcal{D}},+,\| \cdot \|$) is a complete, normed vector space. Moreover, ($\Vf {\mathcal{D}},+, \cdot$) is an associative and commutative $\RR$-algebra with respect to pointwise multiplication. Finally, by Lemma \ref{submult}, the norm is also submultiplicative.
\end{Proof}

\subsection{Further Properties}

In \cite{Blu89} the maximal ideal space of the Banach algebra of functions of bounded Hardy Krause variation was determined.
This result is based on two key observations, \cite[Lemma 6]{Blu89} and \cite[Proposition 3]{Blu89}. In the following we briefly discuss the difficulties in generalising these results to our new notion.\\

We start by recalling the results of \cite{Blu89}. To state the first one, we introduce the function 
\begin{equation*}
\sgn(x)=\begin{cases} 1 & \text{ if } x>0,\\
0 &\text{ if } x=0,\\
-1 &\text{ if } x<0.\\
\end{cases}
\end{equation*}

\begin{lem}[Bl\"{u}mlinger, \cite{Blu89}] \label{blulem6}
Let $f:[0,1]^d \to \R$ be of bounded Hardy--Krause variation. Let $(\mathbf{x}^{(n)})_{n\geq 1}$ be a sequence in $[0,1]^d$ converging to some $\mathbf{x}$ and having the property that $\sgn(x_i^{(n)}-x_i)$ depends only on $i$, but not on $n$. Then $\lim_{n\to \infty} f(\mathbf{x}^{(n)})$ exists.  
\end{lem}

This Lemma was used to prove the following classification of the maximal ideals in $\mathcal{HK}$.

\begin{thm}[Bl\"{u}mlinger, \cite{Blu89}] \label{maxideal}
The maximal ideals in $\mathcal{HK}$ are in one-to-one correspondance with the set of pairs $(\mathbf{x}, {\eps} )\in [0,1]^d \times \{-1,0,1\}^d$ satisfying $0\leq \eps_i$ if $x_i = 0$ and $\eps_i \leq 0$ if $x_i=1$.  The pair $(\mathbf{x},\eps)$ corresponds to the ideal consisting of functions $f\in \mathcal{HK}$ such that $\lim_{\mathbf{x}^{(n)} \to \mathbf{x}} f(\mathbf{x}^{(n)}) = 0$ for any sequence $(\mathbf{x}^{(n)})_{n\geq 0}$ converging to $\mathbf{x}$ and having $\sgn(x_i^{(n)} - x_i) =\eps_i$.
\end{thm}

Lemma \ref{blulem6} does not hold in general for functions of bounded $\mathcal{D}$-variation when $\mathcal{D}\neq \mathcal{R}^*$, as the following example shows.

\begin{ex}
The function $f:[0,1]^2 \to \R$ given by $f(\mathbf{x})= \mathds{1}_{x_1>x_2}$ has bounded $\mathcal{K}$-variation. However, the sequence $\mathbf{x}^{(n)}=(1/2 + 3/n, 1/2 + (2 +(-1)^n)/n )$ converges to $\mathbf{x}=(1/2,1/2)$ and $\sgn(x_i^{(n)} - x_i) = 1$ for all $n$, but $f(\mathbf{x}^{(n)})$ alternates between $0$ and $1$, so $\lim_{n\to \infty} f(\mathbf{x}^{(n)})$ does not exist.
\end{ex} 

Instead, the following weaker version of Lemma \ref{blulem6} holds for $\mathcal{K}$-variation.  For $\mathbf{x}\in [0,1]^d$ and a unit vector $\mathbf{u}\in \Ss^{d-1}$, we define $L_{\mathbf{x},\mathbf{u}}=\{\mathbf{x}+t\mathbf{u}\mid t>0\}$ to be the open half line starting from $\mathbf{x}$ and spanned by $\mathbf{u}$.

\begin{prop}
Suppose $f:[0,1]^d \to \R$ has bounded $\mathcal{K}$-variation. Let $\mathbf{x}\in [0,1]^d$ and $\mathbf{u}\in \Ss^{d-1}$ and let $\mathbf{x}^{(n)}$ be a sequence contained in $L_{\mathbf{x},\mathbf{u}}$  converging to $\mathbf{x}$. Then $\lim_{n\to \infty} f(\mathbf{x}^{(n)})$ exists. 
\end{prop}

\begin{Proof}
Let $K\subseteq [0,1]^d$ be a convex set. Then there is an $\eps >0$ such that either $B_\eps(\mathbf{x})\cap L_{\mathbf{x},\mathbf{u}} \subseteq K$ or $B_\eps(\mathbf{x})\cap L_{\mathbf{x},\mathbf{u}}\subseteq \R^d \backslash K $.  This is obvious if $K\cap L_{\mathbf{x},\mathbf{u}}=\emptyset$. Otherwise, $K\cap L_{\mathbf{x},\mathbf{u}}$ contains a point $\mathbf{y}$. The open line segment between $\mathbf{x}$ and $\mathbf{y}$ is either contained in $K\cap L_{\mathbf{x},\mathbf{u}}$  or it contains a point $\mathbf{y}'=\lambda \mathbf{x} + (1-\lambda )\mathbf{y} $ with $\lambda \in (0,1)$ that does not belong to $K$. By convexity, the line segment between $\mathbf{x}$ and $\mathbf{y}'$ cannot contain any point from $K$. This yields the proposition for $f=\mathds{1}_K$, and hence also when $f$ is a simple function.\\

Let $\mathbf{x}^{(n)}$ be a sequence in $L_{\mathbf{x},\mathbf{u}}$ converging to $\mathbf{x}$ and let $f$ be a general function of bounded $\mathcal{K}$-variation. We must show that $f(\mathbf{x}^{(n)})$ is a Cauchy sequence. Let $\eps>0 $ be given and choose a simple function $g$ with $|f-g|_\infty \leq \eps/3$. Then $|f(\mathbf{x}^{(n)})-f(\mathbf{x}^{(m)})| \leq |g(\mathbf{x}^{(n)})-g(\mathbf{x}^{(m)})| + 2\eps / 3 $. Choose $N $ such that $|g(\mathbf{x}^{(n)})-g(\mathbf{x}^{(m)})|\leq \eps /3$ for all $m,n \geq N$. Then $|f(\mathbf{x}^{(n)})-f(\mathbf{x}^{(m)})| \leq \eps$ for all $m,n \geq N$.
\end{Proof}

The obvious generalization of Theorem \ref{maxideal} would be that the maximal ideals consist of functions whose limits along a fixed line segment $L_{\mathbf{x},\mathbf{u}}$ vanish. However, we have not been able to show this. An important ingredient in the proof of Theorem \ref{maxideal} is Proposition 3 of \cite{Blu89}:

\begin{prop}[Bl\"{u}mlinger, \cite{Blu89}] 
If $f\in \mathcal{HK}$  and there is a $\delta >0$ such that $|f|\geq \delta$, then $1/f \in \mathcal{HK}$.
\end{prop}

While it is relatively simple to see this from the definition of Hardy--Krause variation, it is not obvious to the authors whether a similar statement can be shown for functions of bounded $\mathcal{K}$-variation.  It is true, however, for functions of bounded generalised Harman variation; see  \cite[Corollary 3.17]{PauSva14}. 

\section*{Acknowlegdements}
The first author is supported by a Schr\"odinger scholarship of the Austrian Science Fund (FWF), and by FWF project I 1751-N26. The first and fourth author are supported by FWF projects F 5507 and 5509, which are parts of the Special Research Program \emph{Quasi-Monte Carlo Methods: Theory and Applications}.
The third author is supported by the Centre for Stochastic Geometry and Advanced Bioimaging, funded by the Villum Foundation. \\
Finally, we would especially like to thank the referee for a very careful study of our manuscript that helped to remove various inaccuracies in an earlier version of this paper. 


\end{document}